\documentclass[preprint]{imsart}

\RequirePackage[OT1]{fontenc}
\RequirePackage{amsthm,amsmath}
\RequirePackage[numbers]{natbib}
\RequirePackage[colorlinks,citecolor=blue,urlcolor=blue]{hyperref}

\usepackage{amsmath,amsthm,amssymb,amsfonts}
\usepackage{wasysym,oldgerm,mathrsfs,cancel}
\usepackage{enumerate,stmaryrd}
\usepackage{dsfont,graphicx,yfonts,color,epic,eepic}
\usepackage{xfrac}
\usepackage{makeidx,epstopdf}

\usepackage{mathtools}
\mathtoolsset{showonlyrefs}
\numberwithin{equation}{section}

\pubyear{2019}
\volume{0}
\issue{0}

\startlocaldefs
\numberwithin{equation}{section}

\usepackage{thmtools}

\theoremstyle{plain}
\newtheorem{thm}{Theorem}[section]

\declaretheorem[style=definition,sibling=thm]{lemma}

\newenvironment{Proof}{\noindent\textbf{Proof.\ }}{\hspace*{\fill}$\Box$\medskip}

\newcommand{\Var}{\mathds{V}\mathrm{ar}}

\newcommand{\PR}{\mathds{P}}
\newcommand{\EW}{\mathds{E}}

\newcommand{\dx}{\mathrm{d}}

\newcommand{\deno}{A}

\DeclareMathOperator{\argmax}{argmax}

\endlocaldefs

\usepackage{sectsty}
\sectionfont{\fontsize{18}{16}\bfseries\scshape\selectfont}
\subsectionfont{\fontsize{14}{12}\bfseries\scshape\selectfont}
\subsubsectionfont{\fontsize{12}{12}\bfseries\scshape\selectfont}
\paragraphfont{\fontsize{12}{10}\bfseries\scshape\selectfont}

\begin{document}

\begin{frontmatter}
\title{Minimax $\boldsymbol{L_2}$-Separation Rate in Testing the Sobolev-type Regularity of a Function
}
\runtitle{Testing the Regularity of a Function}

\begin{aug}
\author{\fnms{Maurilio} \snm{Gutzeit}
	\ead[label=e3]{gutzeit@ovgu.de}}

\address{OvGU Magdeburg, Institut f{\"u}r Mathematische Stochastik\\
	Universit{\"a}tsplatz 2, 39106 Magdeburg, Germany\\
\printead{e3}}

\runauthor{M. Gutzeit}

\affiliation{OvGU Magdeburg}

\end{aug}

\begin{abstract}
In this paper we study the problem of testing if an $L_2-$function $f$ belonging to a certain $l_2$-Sobolev-ball $B_t(R)$ of radius $R>0$ with smoothness level $t>0$ indeed exhibits a higher smoothness level $s>t$, that is, belongs to $B_s(R)$. We assume that only a perturbed version of $f$ is available, where the noise is governed by a standard Brownian motion scaled by $\frac{1}{\sqrt{n}}$. More precisely, considering a testing problem of the form
$$H_0:~f\in B_s(R)~~\mathrm{vs.}~~H_1:~f\in B_t(R),~\inf_{h\in B_s}\Vert f-h\Vert_{L_2}>\rho$$
for some $\rho>0$, we approach the task of identifying the smallest value for $\rho$, denoted $\rho^\ast$, enabling the existence of a test $\varphi$ with small error probability in a minimax sense. By deriving lower and upper bounds on $\rho^\ast$, we expose its precise dependence on $n$:
$$\rho^\ast\sim n^{-\frac{t}{2t+1/2}}.$$ As a remarkable aspect of this composite-composite testing problem, it turns out that the rate does not depend on $s$ and is equal to the rate in signal-detection, i.e. the case of a simple null hypothesis. 
\end{abstract}

\begin{keyword}[class=MSC]
\kwd{62G10}
\end{keyword}

\begin{keyword}
	\kwd{minimax hypothesis testing}
	\kwd{nonasymptotic minimax separation rate}
	\kwd{Gaussian white noise}
	\kwd{Sobolev ball}
	\kwd{smoothness}
\end{keyword}
%
\setcounter{tocdepth}{2}
\tableofcontents
\end{frontmatter}

\section{Introduction}\label{s:forProSta}
Let $n\in\mathbb{N}^\ast=\mathbb{N}\backslash\{0\}$, $f$ a fixed unknown element of
$$L_2:=L_2([0,1])=\left\{g:[0,1]\rightarrow\mathbb{R}~;~\int_0^1 g(x)^2~\mathrm{d}\lambda(x)<\infty\right\}$$
and $(B(x))_{x\in[0,1]}$ a standard Brownian motion. Suppose we observe the Gaussian process $(Y(x))_{x\in[0,1]}$ determined by the stochastic differential equation
\begin{equation}\dx Y(x)=f(x)\dx x+\frac{1}{\sqrt{n}}\dx B(x),~~x\in[0,1].\label{s:SDE}\end{equation}
The resulting probability measure, expectation and variance given $f$ will be written $\PR_f$, $\EW_f$ and $\Var_f$, respectively. Depending on the context and if there is no risk of confusion we may drop the index $f$ or write another index, for instance in the context of lower bounds (section \ref{s:LB}).\\
\paragraph{Testing problem}~\\
We now fix $s>t>0$ and $R,\rho>0$. For any $r>0$, we denote by $B_r(R)$ the $l_2$-Sobolev-ball of radius $R$ of functions on $[0,1]$ with regularity at least $r$ -- see section \ref{s:setting} for a precise definition. Based on that, let
$$\widetilde{B}_{s,t}(R,\rho):=\left\{g\in B_t(R)~;~\inf_{h\in B_s(R)}\Vert g-h\Vert_{L_2}>\rho\right\}.$$
Hence, if we interpret $s$ and $t$ as degrees of smoothness, $\widetilde{B}_{s,t}(R,\rho)$ is the set of functions with smoothness level at least $t$ which are separated from the class $B_s(R)$ with stronger smoothness $s$ by $\rho$ in $L_2$-sense. Now, the testing problem of interest is
\begin{equation}H_0:~f\in B_s(R)~~\mathrm{vs.}~~H_1:~f\in \widetilde{B}_{s,t}(R,\rho).\label{s:TestingProblem}\end{equation}
More specifically, given $\eta\in(0,1)$, we aim at finding the magnitude in terms of $n$ of the smallest separation distance $\rho^\ast(\eta)=\rho^\ast(n,s,t,\eta)$ which enables the existence of a test $\varphi$ of level $\eta$ in a minimax sense, i.e. of
\begin{equation}\rho^\ast(\eta)=\inf\left\{\rho>0~;~\exists~\mathrm{test}~\varphi:~\sup_{f\in B_s(R)}\PR_f(\varphi=1)+\sup_{f\in \widetilde{B}_{s,t}(R,\rho)}\PR_f(\varphi=0)\leq\eta\right\}.\label{s:testprob}\end{equation}
\paragraph{Related questions and literature}~\\
There are in essence two lines of work with questions or ideas closely related to the present paper.\\~\\
Firstly, considering the simpler null hypothesis $H_0:~f\equiv 0$ puts us in the so-called signal-detection setting which has already been studied, see for instance the series of seminal papers \cite{Ing93} as well as \cite{Ing02,Lep99}.  or \cite{Com13} for a more recent treatment or \cite{Spo96} for the question of adaptivity to .  In that context, the order of $\rho^\ast$ with respect to $n$ is shown to be
$$n^{-\frac{t}{2t+1/2}}.$$
Moreover, the question of adaptivity to e.g. $t$ is considered in \cite{Spo96} and \cite{Bar02} covers signal detection for Besov balls in a Gaussian sequence setting.\\~\\
Secondly, another closely related task is the construction of (adaptive and honest) confidence regions for $f$. In \cite{Bul13}, the authors study such sets in terms of $L_2$-separation, but rather than the observation $(Y(x))_{x\in[0,1]}$ they use a Gaussian sequence model. However, due to the asymptotic equivalence of these models in the sense of Le Cam (see \cite{Cam12}), it is possible to derive from their arguments that for our problem \eqref{s:TestingProblem},
\begin{equation}n^{-\frac{t}{2t+1/2}}\lesssim\rho^\ast(\eta)\lesssim\max\left(n^{-\frac{s}{2s+1}},n^{-\frac{t}{2t+1/2}}\right).\label{s:resNickl}\end{equation}
While the resulting gap in the case $s<2t$ is not essential in the confidence region setting (see also \cite{Cai06} and \cite{Jud03}), it is quite important from a testing perspective as it raises the question how the complexity of the null hypothesis influences the separation rate.\\
Further relevant literature on confidence sets and adaptivity would be \cite{Car13} as well as \cite{Car18} (matrix completion), \cite{Nic13} (linear regression) and \cite{Gin10} (density estimation). Moreover, interestingly, literature on the frequentist coverage of Bayesian credible sets reveals conditions (``polished tail'', ``self-similarity'' or also ``excessive bias'') which enable deriving adaptive or honest confidence sets from adaptive Bayesian credible sets - see for instance \cite{Sza15} (adaptive confidence sets in Gaussian sequence model with Sobolev-type regularity), \cite{Rou16} (honest confidence sets in rather abstract framework) or also \cite{Bel17} (white-noise model, e.g. adaptive minimax results for the setting from \cite{Sza15} under ``excessive bias'').\\~\\
Now, the article \cite{Car15} is by far the closest previous work to the present paper. Indeed, the author studies the same problem with another choice of Sobolev-ball, namely the $(r,\infty)$-Sobolev-balls $B_{r,\infty}(R)$.
In this context, $\rho^\ast(\eta)$ is proved to be of magnitude
$$n^{-\frac{t}{2t+1/2}}.$$
Note that this quantity is equal to the rate in the signal-detection case and hence in particular does not depend on $s$. This makes the issue of the gap in \eqref{s:resNickl} even more interesting and, from a technical perspective, it is rather striking given that moving from a simple to the composite null hypothesis is a significant step. On top of that, there are settings where the separation rate strongly depends on the shape of the null hypothesis, see e.g. \cite{Bla17} and \cite{Jud02} or also \cite{Cai11}.\\
To the best of our knowledge, the case of \cite{Car15} is the only one for which the minimax $L_2$-separation rate is known and our main contribution is to extend that result to the $(r,2)$-Sobolev-space. While our lower bound (Theorem \ref{s:LB} in section \ref{s:MaiRes}) is essentially a corollary of the corresponding result \cite[Theorem 3.2]{Car15}, the upper bound (Theorem \ref{s:UB} in section \ref{s:MaiRes}) cannot be established through a simple application of \cite[Theorem 3.1]{Car15}. As $B_r(R)\subseteq B_{r,\infty}(R)$, this might be surprising at first sight: Indeed, the test from \cite{Car15} would perform well in the present setting in terms of type-I-error. However, ensuring sufficient power is significantly more difficult when considering $l_2$-Sobolev-balls, see \ref{s:B_infty} for an explicit example.

\section{Setting}\label{s:setting}
In this section, we describe how the relevant Sobolev balls and the observed Gaussian process will be represented throughout the paper.
\paragraph{Wavelet transform and associated Sobolev ball}~\\
Throughout the paper, we make heavy use of a wavelet decomposition of $f$. As is well-known, we can define a scalar product and associated norm on $L_2$ by
$$<g,h>:=\int_0^1 g(x)h(x)~\dx \lambda(x)~~\mathrm{with}~~\Vert g\Vert_{L_2}:=\sqrt{<g,g>},~~g,h\in L_2.$$
There are many orthogonal wavelet bases of $L_2$ with respect to $<\cdot,\cdot>$. A suitable choice for our purposes is a basis developed in \cite{Coh93b} that can be written as
\begin{equation}\mathcal{W}=\bigcup_{j=2}^\infty\{\psi_{j,k}:~k\in\{1,2,\ldots,2^j\}\},\label{eq:wavelet_basis}\end{equation}
i.e. it is tailored such that there are exactly $2^j$ basis functions at resolution $j\geq 2$. Clearly, the coefficients of $g\in L_2$ with respect to $\mathcal{W}$ are given by
$$<g,\psi_{j,k}>=\int_0^1 g(x)\psi_{j,k}(x)~\dx x,~~j\geq 2,k\in\{1,2,\ldots,2^j\}.$$
and yield the representation
\begin{equation}g=\sum_{j=2}^\infty\sum_{k=1}^{2^j}<g,\psi_{j,k}>\psi_{j,k}.\label{s:WaveletCoef}\end{equation}
Let $r>0$. By virtue of isometry properties discussed for instance in \cite{Tri92} and \cite{Nic16}, we may now define a functional $(r,2)$-Sobolev-ball of radius $R$ solely through the wavelet coefficients of its elements, based on the basis from \eqref{eq:wavelet_basis}:
\begin{equation}B_{r}(R):=\left\{g\in L_2~;~\sum_{j=2}^\infty 4^{jr}\sum_{k=1}^{2^j}<g,\psi_{j,k}>^2~\leq~ R^2\right\}\label{s:SobBall}\end{equation}
with associated $(r,2)$-Sobolev-norm
$$\Vert g\Vert_{\mathcal{B}_r}:=\sqrt{\sum_{j=2}^\infty 4^{jr}\sum_{k=1}^{2^j}<g,\psi_{j,k}>^2},~~g\in L_2$$
or also, as mentioned at the end of the previous section,
$$B_{r,\infty}(R):=\left\{g\in L_2~;~\sup_{j\geq 2}~ 4^{jr}\sum_{k=1}^{2^j}<g,\psi_{j,k}>^2~\leq~ R^2\right\}.$$
\paragraph{Discrete observation scheme based on the wavelet basis}~\\
Let
$$\mathcal{I}=\{(j,k)\in\mathbb{N}^2~|~j\geq 2,k\leq 2^j\}.$$
Motivated by \eqref{s:SobBall}, for each $(j,k)\in\mathcal{I}$ we consider
$$a_{j,k}:=<f,\psi_{j,k}>$$
so that
$$f=\sum_{j=2}^\infty\sum_{k=1}^{2^j}a_{j,k}\psi_{j,k}.$$
The natural corresponding estimators read
\begin{equation}\widehat{a}_{j,k}:=<\dx Y,\psi_{j,k}>,~~\widehat{f}=\sum_{j=2}^\infty\sum_{k=1}^{2^j}\widehat{a}_{j,k}\psi_{j,k}.\label{s:estimators}\end{equation}
By construction and due to the orthonormality of $\mathcal{W}$, we know that the family $(\widehat{a}_{i,j})_{(j,k)\in\mathcal{I}}$ is independent with
$$\widehat{a}_{j,k}\sim\mathcal{N}\left(a_{j,k},\tfrac{1}{n}\right).$$
Clearly, observing this family is equivalent to observing the original process $(Y(x))_{x\in[0,1]}$.
\section{Main results}\label{s:MaiRes}
In this section, we state and discuss our main results, that is upper and lower bounds on $\rho^\ast(\eta)$. We also provide a high-level description of the strategy and ideas included in the upper bound proof, which is our main contribution.
\subsection{Upper Bound}
\paragraph{The test}~\\
Note that $\widehat{f}$ from \eqref{s:estimators} is not a useful estimator as it exhibits infinite variance. Therefore, we need to carefully impose a restriction of the form $j\leq J$ for some fixed $J\in\mathbb{N}$, $J\geq 2$. Actually, section \ref{s:UBchapter} is primarily concerned with obtaining an upper bound on $\rho_J^\ast(\eta)$ for the reduced, finite-dimensional problem
$$H_0':~\underbrace{\left\Vert \sum_{j=2}^J\sum_{k=1}^{2^j}a_{j,k}\psi_{j,k}\right\Vert_{\mathcal{B}_s}}_{:=S_J}\leq R~~~\mathrm{vs.}~~~H_1':~\inf_{h\in B_s(R)}\left\Vert \sum_{j=2}^J\sum_{k=1}^{2^j}a_{j,k}\psi_{j,k}-h\right\Vert_{L_2}>\rho_J,$$
where $\rho_J$ and $\rho_J^\ast(\eta)$ are analogous in definition and relation to their counterparts in \eqref{s:TestingProblem} and \eqref{s:testprob}. In fact, finding a sufficient separation distance $\rho_J\geq\rho^\ast_J(\eta)$ here is the central and most involved part of the paper.\\
As we illustrate in section \ref{s:B_infty}, it turns out that a test based on estimating $S_J^2$ only cannot perform well enough under the targeted separation distance of order $n^{-t/(2t+1/2)}$ due to the strong variance at high levels, so that more flexibility is necessary: In Lemma \ref{s:Jpart2}, we analyse the smallest level $j^\ast$ such that $S_{j^\ast}$ considerably exceeds $R$ (such an index must exist under $H_1'$) and it turns out that this is detectable through the estimator $\Vert P_2^{j^{\ast}}\widehat{f}\Vert_{\mathcal{B}_{s}}^2$ (section \ref{s:conclusion}, second paragraph). Hence, we propose a test which evaluates the individual accumulated (squared) Sobolev-norms of the projections up until level $J$ and rejects the null hypothesis whenever one of these norms is too large.\\~\\
In particular, we define for $j^\ast\in\{2,3,\ldots,J\}$
\begin{align*}
\alpha_{j^\ast}&=\eta\frac{1-2^{-1/5}}{4}2^{(j^\ast-J)/5},\\
T_{j^\ast,\alpha_{j^\ast}}&=\Vert P_2^{j^\ast}\widehat{f}\Vert_{\mathcal{B}_s}^2-A_{j^\ast}-\frac{2}{\sqrt{\alpha_{j^\ast}}}\cdot \frac{\sqrt{j^\ast-1}}{\sqrt{n}}\sqrt{\max_{2\leq j\leq j^\ast}\left|16^{js}\left(\Vert P_j\widehat{f}\Vert_{L_2}^2-\frac{2^j}{n}\right)\right|},\\
\beta_{j^\ast}&=\eta\frac{1-2^{-1/2}}{2}2^{-j^\ast/2},\\
	C_{\beta_{j^\ast}}&=\sqrt{\frac{2}{\beta_{j^\ast}}},\\
D_{j^\ast,\beta_{j^\ast}}&=\frac{4^{j^\ast s}}{\sqrt{n}}\left(\sqrt{2}C_{\beta_{j^\ast}}+2^{j^\ast/4}\sqrt{C_{\beta_{j^\ast}}}\right),\\
\tau_{j^\ast,\alpha_{j^\ast}}&=R^2+\frac{2}{\sqrt{\alpha_{j^\ast}}}\left( \frac{\sqrt{j^\ast-1}}{\sqrt{n}}D_{j^\ast,\beta_{j^\ast}}+4^{j^\ast s}\frac{2^{j^\ast/2}}{n}\right)
\end{align*}
and finally the test
\begin{equation}\varphi=1-\prod_{j^{\ast}=2}^{J}\mathds{1}_{\{T_{j^\ast,\alpha_{j^\ast}}\leq\tau_{j^\ast,\alpha_{j^\ast}}\}}.\label{test}\end{equation}
In principle, the conditions $T_{j^\ast,\alpha_{j^\ast}}\leq\tau_{j^\ast,\alpha_{j^\ast}}$ are based on applying Chebyshev's inequality to the estimators $\Vert P_2^{j^{\ast}}\widehat{f}\Vert_{B_{s}}^2$ with a bias-correction term $A_{j^\ast}$ (Lemma \ref{s:concLemma} below). Now, since the variance of $\Vert P_2^{j^{\ast}}\widehat{f}\Vert_{B_{s}}^2$ depends on $f$, it needs to be estimated, which manifests itself especially in the last part of $T_{j^\ast,\tau_{j^\ast}}$.\\~\\
The choice of $J$ is then governed by reaching a trade-off between the resulting upper bound on $\rho_J^\ast(\eta)$ and the error incurred by ignoring the resolutions beyond $J$ - it is the index where they are both of order $n^{-\frac{t}{2t+1/2}}$,
$$J=\left\lfloor \frac{1}{2t+1/2}\frac{\ln(n)}{\ln(2)}\right\rfloor.$$\\~\\
In terms of technical ingredients, all these considerations are remarkable in that they solely rely on elementary computations based on the Sobolev-balls' geometry and classical properties of the $\chi^2-$distribution.\\~\\
Our main result reads as follows:
\begin{thm}\label{s:UB}
	Let $\eta\in(0,1)$. Whenever
	$$\rho\geq \left(\frac{1346}{\sqrt{\eta}}+\frac{R}{1-2^{-t}}\right)n^{-\frac{t}{2t+1/2}}\label{s:rhofinal},$$
	the test $\varphi$ from \eqref{test} fulfils
	$$\sup_{f\in B_s(R)}\PR_f(\varphi=1)+\sup_{f\in \widetilde{B}_{s,t}(R,\rho)}\PR_f(\varphi=0)\leq\eta.$$
	Hence,
	$$\rho^\ast(\eta)\leq \left(\frac{1346}{\sqrt{\eta}}+\frac{R}{1-2^{-t}}\right)n^{-\frac{t}{2t+1/2}},~~\mathrm{i.e.}~~\rho^\ast(\eta)\lesssim n^{-\frac{t}{2t+1/2}}.$$
\end{thm}

\subsection{Remark on the relation to \cite{Car15}}\label{s:B_infty}
In order to clarify the distinction between the previous work \cite{Car15} with $H_0:~f\in B_{s,\infty}(R)$ and the present paper, we consider two rather specific examples.

\paragraph{Testing the resolutions separately does not suffice}~\\
First of all, note that $B_{s,\infty}(R)$ is very large compared to $B_{s}(R)$, which ensures that, as mentioned above, the test $\Psi$ from \cite{Car15} performs well under the null hypothesis $H_0$ of the present paper. However, this geometric imbalance is so strong that often for one and the same function, we would like one test to reject the null hypothesis and the other test to not reject it:\\~\\
Consider a simple extreme case where
$$\Vert P_j f\Vert_{L_2} = \frac{R}{2^{js}},~~j\in\{2,3,\ldots\}.$$
Then clearly we have
\begin{eqnarray*}\Vert f\Vert_{\mathcal{B}_{s,\infty}} &=& R,\\
	\Vert f\Vert_{\mathcal{B}_{s}} &=& \sqrt{\sum_{j=2}^\infty R^2} = \infty,\\
	\Vert f\Vert_{\mathcal{B}_{t}} &=& R\sqrt{\sum_{j=2}^\infty 4^{j(t-s)}}.
\end{eqnarray*}
It can be assured that $f\in B_t(R)$ through the condition $t< s-\log_4\left(\frac{2}{\sqrt{5}-1}\right)$, so that clearly we have found a case
where
$$f\in B_{s,\infty}(R),~~f\notin B_{s}(R),~~f\in B_{t}(R),$$
i.e. both the null hypothesis of \cite{Car15} and our alternative hypothesis are met. The test from \cite{Car15} based on separately evaluating the individual levels will clearly not reject our null hypothesis with high probability. On the other hand, in order to check the new test's performance, let us invoke Theorem \ref{s:UB}: By construction, for any $h\in B_{s}(R)$, there is a sequence $(a_j)_{j\in\{2,3,\ldots\}}$ in $[0,1]$ such that
$$\Vert P_j h\Vert_{L_2} = a_j\frac{R}{2^{js}}~\mbox{and}~\sum_{j=2}^\infty a_j^2 = 1.$$ Then we have
\begin{eqnarray*}\Vert f-h\Vert_{L_2}^2 &=& \sum_{j=2}^\infty \Vert P_jf-P_jh\Vert_{L_2}^2\\
	&\geq& \sum_{j=2}^\infty\left(\Vert P_jf\Vert_{L_2}-\Vert P_jh\Vert_{L_2}\right)^2\\
	&=& \sum_{j=2}^\infty \frac{R^2}{4^{js}}\left(1-a_j\right)^2\\
	&´\geq& R^2 \frac{3-2\sqrt{2}}{4^{3s}},\end{eqnarray*}
where the last bound can be derived from the observation that necessarily $a_2\leq \frac{1}{\sqrt{2}}$ or $a_3\leq\frac{1}{\sqrt{2}}$.
As this holds for any $h\in B_{s}(R)$, in particular we have
$$\inf_{h\in B_{s}(R)}\Vert f-h\Vert_{L_2}~\gtrsim~R~\gtrsim~ n^{-\frac{t}{2t+1/2}}$$
for appropriate $n\in\mathbb{N}$ so that the new test detects that $f\notin B_{s}(R)$ with high probability.
\paragraph{Estimating only $\boldsymbol{\Vert P_2^J f\Vert_{\mathcal{B}_{s,2}}}$ does not suffice}~\\
The strategy of only estimating $\Vert P_2^J f\Vert_{\mathcal{B}_{s}}$ is too optimistic in the present setting:\\~\\
Consider a case where for some $a>1$
$$\Vert P_2f\Vert_{L_2}^2 = a^2\frac{R^2}{4^{2s}},~~\Vert P_Jf\Vert_{L_2}^2 = \frac{R^2}{4^{Js}},~~\Vert P_jf\Vert_{L_2} = 0~~\mbox{else}.$$
Then on the one hand,
\begin{eqnarray*}
\inf_{h\in B_{s}(R)}\Vert f-h\Vert_{L_2}^2 &=& \inf_{h\in B_{s}(R)}\left(\Vert P_2(f-h)\Vert_{L_2}^2 + \Vert P_J(f-h)\Vert_{L_2}^2\right)\\
&\geq& \left(\Vert P_2f\Vert_{L_2} - \sup_{h\in B_{s}(R)}\Vert P_2f\Vert_{L_2}\right)^2\\
&=&\frac{(a-1)^2}{4^{2s}}R^2,
\end{eqnarray*}
which, again, exceeds our (squared) upper bound for appropriate $n$ or $a$ so that in principle, it is possible to detect $f\notin B_{s}(R)$ in the sense of Theorem \ref{s:UB}.\\
Note that we can see this without using information on more than the first level. This is an important observation with regards to the construction of our test.\\
Furthermore, we have
$$\inf_{h\in B_{s}(R)}\Vert f-h\Vert_{L_2}^2 \leq \Vert f\Vert_{L^2}^2 = \left(\frac{a^2}{4^{2s}}+\frac{1}{4^{Js}}\right)R^2.$$
On the other hand, as we show in Lemma \ref{s:concLemma} below, in this special case the cost in terms of standard deviation of including the estimate $\Vert P_J\widehat{f}\Vert_{\mathcal{B}_{s}}^2$ would be
\begin{eqnarray*}\sqrt{\Var\left[\Vert P_J\widehat{f}\Vert_{\mathcal{B}_{s}}^2\right]}&\gtrsim&4^{Jt}\frac{2^{J/2}}{n}\\
&\gtrsim& n^{0} = 1\end{eqnarray*}
(absolute constand). For large enough $s$ and/or small enough $a$, this standard deviation exceeds the (squared) distance to be detected - hence, a test based on level $J$ is unlikely to correctly reject the null hypothesis. The test we propose copes with such a situation through analysing multiple accumulated estimates and would  have detected $f\in B_{s}(R)$ at the first level already with high probability.
\subsection{Lower Bound}
Using the same choice for $J$ as indicated above, a lower bound on $\rho^\ast(\eta)$ of the same order can be derived through studying the statistical distance between specific distributions agreeing with $H_0$ and $H_1$ respectively.
\begin{thm}\label{s:LB}
	Let $\eta\in(0,1)$. There are $C_\eta>0$ and $N_\eta\in\mathbb{N}$ such that whenever $n\geq N_\eta$ and
	$$\rho\leq C_{\eta}n^{-\frac{t}{2t+1/2}},$$
	for any test $\varphi$ it holds that
	$$\sup_{f\in B_s(R)}\PR_f(\varphi=1)+\sup_{f\in \widetilde{B}_{s,t}(R,\rho)}\PR_f(\varphi=0)>\eta.$$
	Hence,
	$$\rho^\ast(\eta)\geq C_{\eta}n^{-\frac{t}{2t+1/2}},~~\mathrm{i.e.}~~\rho^\ast(\eta)\gtrsim n^{-\frac{t}{2t+1/2}}.$$
	In particular, one may choose
	$$C_\eta:=\frac{R}{2}\min\left\{1,\frac{\sqrt{\ln(1+4(1-\eta)^2)}}{2^t 16R}\right\},~~N_\eta:=\left\lceil \left(R\frac{2^{s-t}}{C_\eta}\right)^{\frac{2t+1/2}{s-t}}\right\rceil.$$
\end{thm}
\noindent Note that, as mentioned in the introduction, Theorems \ref{s:UB} and \ref{s:LB} in conjunction reveal the minimax separation rate to be of order
$$\rho^\ast(\eta)\sim n^{-\frac{t}{2t+1/2}},$$
which does not depend on the size of the null hypothesis and is equal to the signal-detection rate. Indeed, in order to obtain the lower bound of Theorem \ref{s:LB}, the fact that $H_0$ is a composite hypothesis need not be used.
\section{Alternative settings}
Before presenting the proofs of our main results, we briefly discuss their possible application in two alternative settings which might also be of interest, see also \cite[Section 3.3]{Car15} and references therein.
\paragraph{Heteroscedastic noise}~\\
As a generalisation of \eqref{s:SDE}, consider the model
\begin{equation}\dx Y(x)=f(x)\dx x+\frac{\sigma(x)}{\sqrt{n}}\dx B(x),~~x\in[0,1],\label{s:hetero}\end{equation}
where $\sigma\in L_2$ is unknown. The proof of Theorem \ref{s:UB} relies heavily on unbiased estimators of $a_{j,k}^2$, $(j,k)\in\mathcal{I}$, and hence on knowledge of the noise coefficient, so that in this generalised version we cannot directly apply our result. However, there is a relatively simple solution under certain conditions: Suppose we have access to two independent realisations $(Y^{(1)}(x))_{x\in[0,1]}$ and $(Y^{(2)}(x))_{x\in[0,1]}$ with noise coefficient, say, $\frac{\sigma(x)}{\sqrt{n/2}}$. Then we can still consider the estimates
$$\widehat{a}^{(i)}_{j,k}~=~<\dx Y^{(i)},\psi_{j,k}>~\sim~\mathcal{N}\left(a_{j,k},2\frac{\Vert \sigma\cdot\psi_{j,k}\Vert_{L_2}^2}{n}\right),~~i\in\{1,2\}$$
and define a new unbiased estimator for $a_{j,k}^2$ based on the simple observation
$$\EW[\widehat{a}^{(1)}_{j,k}\cdot \widehat{a}^{(2)}_{j,k}]=a_{j,k}^2.$$
If in addition we know an upper bound on $\Vert \sigma\Vert_{L_2}$, it turns out that we can state an analogous concentration result as the one for the homoscedastic model (see Lemma \ref{s:concLemma} below) and obtain essentially the same result.
\paragraph{Regression}~\\
Another possible observation scheme for testing the smoothness of $f$ would be collecting $n$ iid samples $(X_i,Y_i)_{i\in\{1,2,\ldots,n\}}$ according to the model
$$Y=f(X)+\frac{\sigma(X)}{\sqrt{n}}\epsilon$$
for $\epsilon\sim\mathcal{N}(0,1)$ and $X$ uniformly distributed on $[0,1]$. This situation is particularly interesting since, as mentioned above, it is asymptotically equivalent to \eqref{s:hetero} in the sense of Le Cam (\cite{Cam12}) We could then arrive at the same situation as in the previous setting by considering
$$\widehat{a}^{(1)}_{j,k}=\frac{2}{n}\sum_{i=1}^{n/2}Y_i\psi_{j,k}(X_i),~~\widehat{a}^{(2)}_{j,k}=\frac{2}{n}\sum_{i=n/2+1}^{n}Y_i\psi_{j,k}(X_i).$$
Note that if $X$ is not uniformly distributed, $\EW[\widehat{a}^{(i)}_{j,k}]=a_{j,k}$ is generally not true and it becomes crucial to guarantee a certain spread of the design points $(X_i)_{i\in\{1,2,\ldots,n\}}$ over $[0,1]$.
\paragraph{Open problems: separation in $\boldsymbol{L}_p$-norm and more general Sobolev-spaces}~\\
We only consider separation in $L_2$-norm, which raises the question if it is possible to generalise the results to separation in $L_p$-norm, $p>0$; the same is true for the previous paper \cite{Car15}. We believe that the strategies of both papers cannot be easily generalised as different values for $p$ result in fundamentally different problems. Indeed, strong differences with varying $p$ already show in the allegedly simple setting of signal detection in the Gaussian vector/sequence model, see \cite[section 3.3.2]{Ing02}). Much more closely related to the present paper, in \cite{Lep99} the authors derive optimal rates for estimating $\Vert f\Vert_{L_p}$ and give very different results and approaches for even versus odd integers $p$. With that said, considering more general Sobolev-balls would seem to produce similar effects as our results heavily rely on estimating the $2$-norm of projections of $f$ (or, in some sense, $\sup$-norm in the previous paper); coping with different parameters here is not trivial as can be seen for instance in the proofs of \cite{Car13}.\\
In summary, such considerations are generally possible and constitute worthwile future work, but they are beyond the scope of the present paper.

\section{Proof of Theorem \ref{s:UB}}\label{s:UBchapter}
\subsection{General preparations}\label{s:GenPrep}
\paragraph{Reduction of the range of resolutions}~\\
Let us make this more clear at this point already: For $j_1,j_2\in\mathbb{N}\cup\{\infty\}$ with $2\leq j_1\leq j_2$ and $g\in L_2$, define the projections
$$P_{j_1}^{j_2}g=\sum_{j=j_1}^{j_2}\sum_{k=1}^{2^j}<g,\psi_{j,k}>\psi_{j,k},~~P_{j_1}:=P_{j_1}^{j_1}.$$
Now observe that since $f\in B_t(R)$, for each $j\in\mathbb{N}$, $j\geq 2$, we have
$$\Vert P_jf\Vert_{L_2}=\frac{\Vert P_j f\Vert_{\mathcal{B}_t}}{2^{jt}}\leq\frac{R}{2^{jt}}$$
and hence
$$\sum_{j=J+1}^{\infty}\Vert P_j f\Vert_{L_2}\leq R\sum_{j=J+1}^\infty (2^{-t})^j=2^{-tJ}\frac{2^{-t}R}{1-2^{-t}}.$$
Using the triangle inequality, this tells us that under the alternative hypothesis
\begin{align}
\rho&<\inf_{h\in B_s(R)}\Vert f-h\Vert_{L_2}\\
&\leq\inf_{h\in B_s(R)}\Vert P_2^J f-h\Vert_{L_2}+\Vert P_{J+1}^\infty f\Vert_{L_2}\\
&\leq\inf_{h\in B_s(R)}\Vert P_2^J f-h\Vert_{L_2}+\sum_{j=J+1}^{\infty}\Vert P_j f\Vert_{L_2},\\
&\leq \inf_{h\in B_s(R)}\Vert P_2^J f-h\Vert_{L_2}+2^{-tJ}\frac{2^{-t}R}{1-2^{-t}}.
\end{align}
Accordingly, under $H_1$ we consider the assumption 
\begin{equation}\rho-2^{-tJ}\frac{2^{-t}R}{1-2^{-t}}=:\rho_J<\inf_{h\in B_s(R)}\Vert P_2^J f-h\Vert_{L_2}\label{s:H1alt}\end{equation}
and firstly solve $\eqref{s:testprob}$ for $\rho_J$ in terms of the reduced range $j\in\{2,3,\ldots,J\}$, that is, subsequently, we will primarily study the testing problem
\begin{equation}H_0':~\Vert P_2^J f\Vert_{\mathcal{B}_s}\leq R~~~\mathrm{vs.}~~~H_1':~\inf_{h\in B_s(R)}\Vert P_2^J f-h\Vert_{L_2}>\rho_J.\label{s:newProb}\end{equation}
Finally, $\rho$ will be determined by choosing $J$ such that a reasonable trade-off between the two summands, 
\begin{equation}\rho=\rho_J+2^{-tJ}\frac{2^{-t}R}{1-2^{-t}},\label{s:rhoDef}\end{equation}
is realised.\\ \ \\
Now, more specifically, with $a=1346$, for $j^\ast\in\{2,3,\ldots,J\}=:\mathcal{J}$, let
\begin{equation}\rho_1:=0;~~\rho_{j^\ast}=a\frac{2^{(3j^\ast+2J)/20}}{\sqrt{n}}.\label{s:eq:rhoDef}\end{equation}
Under $H_1'$ it will be technically useful to detect the level $\boldsymbol{j^\ast}\in\mathcal{J}$ at which $\displaystyle\inf_{h\in B_s(R)}\Vert P_2^{j^\ast}f-h\Vert_{\mathcal{B}_s}$ firstly exceeds $\rho_{j^\ast}$ in the sense of Lemma \ref{s:LevelLemma} below. That leads to a multiple test across the set $\mathcal{J}$ finally given in \eqref{s:thetest}.
\paragraph{Decomposition of $\boldsymbol{H_1'}$}
\begin{lemma}\label{s:LevelLemma}Under the alternative hypothesis $H_1'$, we have
	\begin{equation}\exists j^\ast\in\mathcal{J}:~\left\{\begin{array}{l}\displaystyle\inf_{h\in B_s(R)}\Vert P_2^{j^\ast-1}f-h\Vert_{L_2}\leq \rho_{j^\ast-1},\\ \displaystyle\inf_{h\in B_s(R)}\Vert P_2^{j^\ast }f-h\Vert_{L_2}> \rho_{j^\ast}.\end{array}\right.\label{s:Logic}\end{equation}
\end{lemma}
\begin{Proof}		
	By contradiction: Assume that \eqref{s:Logic} is false, i.e.
	$$\forall j^\ast\in\mathcal{J}:~~\underbrace{\inf_{h\in B_s(R)}\Vert P_2^{j^\ast-1}f-h\Vert_{L_2}> \rho_{j^\ast-1}}_{E_{j^\ast}}~\vee~\underbrace{\inf_{h\in B_s(R)}\Vert P_2^{j^\ast }f-h\Vert_{L_2}\leq \rho_{j^\ast}}_{F_{j^\ast}}.$$
	Then clearly $F_J$ is false, so that $E_{J}$ is true. Equivalently, $F_{J-1}$ is false and in turn $E_{J-1}$ must be true. Continued application of this argument leads to the contradiction $$\displaystyle\inf_{h\in B_s(R)}\Vert P_2^{1}f-h\Vert_{L_2}=0> \rho_{1}.$$
\end{Proof}
\paragraph{Concentration of $\boldsymbol{\Vert P_2^{j^\ast} \widehat{f}\Vert_{\mathcal{B}_s}^2}$}
\begin{lemma}\label{s:concLemma}
	Let $j^\ast\in\mathcal{J}$. Then, with
	\begin{equation}
	A_{j^\ast}:=\frac{1}{n}\sum_{j=2}^{j^\ast}(2\cdot 4^s)^j,~B_{j^\ast}:=\frac{2}{n^2}\sum_{j=2}^{j^\ast}(2\cdot 4^{2s})^j,~	V_{j^\ast}=\frac{4}{n}\sum_{j=2}^{j^\ast}4^{2js}\Vert P_jf\Vert_{L_2}^2,
	\label{s:ConcPara}
	\end{equation}
	it holds that
	\begin{equation}\forall\delta\in(0,1):~~
	\PR\left(\left| \Vert P_2^{j^\ast}\widehat{f}\Vert_{\mathcal{B}_s}^2-A_{j^\ast}-\Vert P_2^{j^\ast}f\Vert_{\mathcal{B}_s}^2\right|\geq \sqrt{\frac{1}{\delta}\left(B_{j^\ast} +V_{j^\ast}\right)}\right)\leq \delta.
	\label{s:eq:ConcIn}
	\end{equation}
\end{lemma}
\begin{Proof}
	For $j\in\mathcal{J}$, let
	$$\lambda_j:=n\sum_{k=1}^{2^j}a_{j,k}^2=n\Vert P_jf\Vert_{L_2}^2.$$
	Then, by construction, we know that
	$$n\Vert P_j\widehat{f}\Vert_{L_2}^2=\sum_{k=1}^{2^j}(\sqrt{n}~\widehat{a}_{j,k})^2\sim\chi^2(2^j,\lambda_j),$$
	i.e. a $\chi^2-$distribution with $2^j$ degrees of freedom and non-centrality parameter $\lambda_j$. Classical properties of this distribution now tell us
	\begin{equation}
	\EW\left[\Vert P_j\widehat{f}\Vert_{L_2}^2\right]=\frac{2^j}{n}+\Vert P_j f\Vert_{L_2}^2;~~\Var\left[\Vert P_j\widehat{f}\Vert_{L_2}^2\right]=2\left(\frac{2^j}{n^2}+\frac{2}{n}\Vert P_j f\Vert_{L_2}^2\right).
	\label{s:eq:Chi2Para}
	\end{equation} 
	Since $$\Vert P_2^{j^\ast}f\Vert_{\mathcal{B}_s}^2=\sum_{j=2}^{j^\ast}4^{js}\Vert P_j\widehat{f}\Vert_{L_2}^2,$$ 
	independence in conjunction with \eqref{s:eq:Chi2Para} yields
	\begin{align}
	\EW[\Vert P_2^{j^\ast}\widehat{f}\Vert_{\mathcal{B}_s}^2]&=\sum_{j=2}^{j^\ast}4^{js}\left(\frac{2^j}{n}+\Vert P_j f\Vert_{L_2}^2\right)\\
	&=\frac{1}{n}\sum_{j=2}^{j^\ast}(2\cdot 4^s)^j+\Vert P_2^{j^\ast}f\Vert_{\mathcal{B}_s}^2\\
	&=A_{j^\ast}+\Vert P_2^{j^\ast}f\Vert_{\mathcal{B}_s}^2;\\
	\Var[\Vert P_2^{j^\ast}\widehat{f}\Vert_{\mathcal{B}_s}^2]
	&=\sum_{j=2}^{j^\ast}4^{2js}\left(2\left(\frac{2^j}{n^2}+\frac{2}{n}\Vert P_j f\Vert_{L_2}^2\right)\right)\\
	&=\frac{2}{n^2}\sum_{j=2}^{j^\ast}(2\cdot 4^{2s})^j+\frac{4}{n}\sum_{j=2}^{j^\ast}4^{2js}\Vert P_jf\Vert_{L_2}^2\\
	&=B_{j^\ast}+V_{j^\ast}.
	\end{align}
	We obtain the desired result directly through Chebyshev's inequality: For $\epsilon>0$,
	$$
	\PR\left(\left| \Vert P_2^{j^\ast}\widehat{f}\Vert_{\mathcal{B}_s}^2-A_{j^\ast}-\Vert P_2^{j^\ast}f\Vert_{\mathcal{B}_s}^2\right|\geq \epsilon\right)\leq \frac{B_{j^\ast} +V_{j^\ast}}{\epsilon^2}
	$$
	and hence the claim.
\end{Proof}~\\
%
More specifically, observe that
\begin{align}
B_{j^\ast}&=\frac{2}{n^2}\sum_{j=2}^{j^\ast}(2\cdot 4^{2s})^j\\
&=\frac{2}{n^2}(2\cdot 4^{2s})^2\frac{(2\cdot 4^{2s})^{j^\ast-1}-1}{2\cdot 4^{2s}-1}\\
&\leq \frac{2}{n^2}\frac{2\cdot 4^{2s}}{2\cdot 4^{2s}-1}(2\cdot 4^{2s})^{j^\ast}\\
&\leq \frac{4}{n^2}(2\cdot 4^{2 s})^{j^\ast},
\end{align}
(where we use that for $x\geq 2$, $\frac{x}{x-1}\leq 2$) and hence for $\delta\in(0,1)$
\begin{equation}\sqrt{\frac{B_{j^\ast}}{\delta}}\leq \frac{2}{\sqrt{\delta}}4^{j^\ast s}\frac{2^{j^\ast/2}}{n}.\label{s:eq:B_j}\end{equation}
Furthermore,
\begin{equation}
\sqrt{\frac{V_{j^\ast}}{\delta}}=\frac{2}{\sqrt{\delta}}\cdot \frac{1}{\sqrt{n}}\sqrt{\sum_{j=2}^{j^\ast}4^{2js}\Vert P_jf\Vert_{L_2}^2}\leq \frac{2}{\sqrt{\delta}}\cdot \frac{\sqrt{j^\ast-1}}{\sqrt{n}}\max_{2\leq j\leq j^\ast}2^{js}\Vert P_jf\Vert_{\mathcal{B}_s}.\label{s:VarBound}
\end{equation}
The maximum in the latter computation will play an important role in the sequel. From now on we use the abbreviation
\begin{equation}
M_{j^\ast}:=\max_{2\leq j\leq j^\ast}2^{js}\Vert P_jf\Vert_{\mathcal{B}_s}.
\label{s:MaxDef}
\end{equation}
Plugging these bounds in \eqref{s:eq:ConcIn} leads to
\begin{equation}\PR\left(\left| \Vert P_2^{j^\ast}\widehat{f}\Vert_{\mathcal{B}_s}^2-A_{j^\ast}-\Vert P_2^{j^\ast}f\Vert_{\mathcal{B}_s}^2\right|\geq\frac{2}{\sqrt{\delta}}\cdot \frac{\sqrt{j^\ast-1}}{\sqrt{n}}M_{j^\ast}+\frac{2}{\sqrt{\delta}}4^{j^\ast s}\frac{2^{j^\ast/2}}{n}\right)\leq \delta
\label{s:eq:ConcIn2}
\end{equation}
for any $\delta\in(0,1)$.
\subsection{Preliminary Bounds on $\boldsymbol{\Vert P_2^{j^\ast}f\Vert_{\mathcal{B}_s}}$}
As a next step towards controlling the type-I and type-II errors of our test, we study $\Vert P_2^{j^\ast}f\Vert_{\mathcal{B}_s}$ more closely.\\
On the one hand, under $H_0'$, for any $j^\ast\in\mathcal{J}$ we clearly have $\Vert P_2^{j^\ast}f\Vert_{\mathcal{B}_s}\leq R$.\\
On the other hand, under $H_1'$, we require a lower bound on $\Vert P_2^{j^\ast}f\Vert_{\mathcal{B}_s}$. The following bound is preliminary in the sense that it requires the knowledge of an index $j^\ast\in\mathcal{J}$ with the property from \eqref{s:Logic} and the corresponding $M_{j^\ast}$. The generalisation will be considered in sections \ref{s:estimation} and \ref{s:conclusion}.
\begin{lemma}\label{s:Jpart2}Let $j^\ast\in\mathcal{J}$ be an index with the property
	\begin{equation}\inf_{h\in B_s(R)}\Vert P_2^{j^\ast-1}f-h\Vert_{L_2}\leq\rho_{j^\ast-1},~~\inf_{h\in B_s(R)}\Vert P_2^{j^\ast}f-h\Vert_{L_2}>\rho_{j^\ast}.\label{s:Logic2}\end{equation}
	Then the following assertion holds for $A=11$:
	\begin{equation}
	\Vert P_2^{j^\ast}f\Vert_{\mathcal{B}_s}^2\geq R^2+\frac{1}{2\cdot \deno^2}\rho_{j^\ast}M_{j^\ast}+\frac{1}{2\cdot\deno^2}4^{j^\ast s}\rho_{j^\ast}^2.
	\label{s:LBH1}
	\end{equation} 
\end{lemma}
\begin{Proof}
	%
	Before giving the main arguments, we need a technical preparation and a general (i.e. only depending on $j^\ast$) lower bound on $\Vert P_2^{j^\ast}f\Vert_{\mathcal{B}_s}$:
	\begin{enumerate}
		\item \underline{Proxy minimisation of $\inf_{h\in B_s(R)}\Vert P_2^{j^\ast}f-h\Vert_{L_2}$}\\
		For $\widetilde{j}\in\mathcal{J}$, write $P_{j\neq\widetilde{j}}:=P_{2}^{j^\ast}-P_{\widetilde{j}}$. In the case that $\Vert P_{j\neq\widetilde{j}}f\Vert_{\mathcal{B}_s}\leq R$, we can introduce the function $\widetilde{h}$ through the wavelet coefficients
		$$\begin{array}{rll}b_{j,k}&:=a_{j,k}&\mathrm{for}~(j,k)\in\mathcal{I},j\neq \widetilde{j},\\ b_{\widetilde{j},k}&:=a_{\widetilde{j},k}\cdot\displaystyle\frac{\sqrt{R^2-\Vert P_{j\neq\widetilde{j}}f\Vert_{\mathcal{B}_s}^2}}{\Vert P_{\widetilde{j}}f\Vert_{\mathcal{B}_s}},&\mathrm{for}~k\in\{1,2,\ldots,2^{\widetilde{j}}\}.\end{array}$$
		Then $\widetilde{h}\in B_s(R)$ holds since
		\begin{align}
		\Vert \widetilde{h}\Vert_{\mathcal{B}_s}^2&=\sum_{j=2}^{j^\ast}4^{js}\sum_{k=1}^{2^j} b_{j,k}^2\\
		&=\Vert P_{j\neq\widetilde{j}}f\Vert_{\mathcal{B}_s}^2+\left(\frac{\sqrt{R^2-\Vert P_{j\neq\widetilde{j}}f\Vert_{\mathcal{B}_s}^2}}{\Vert P_{\widetilde{j}}f\Vert_{\mathcal{B}_s}}\right)^2\Vert P_{\widetilde{j}}f\Vert_{\mathcal{B}_s}^2\\
		&=R^2.
		\end{align}
		Hence, by assumption
		$$\rho_{j^\ast}^2< \inf_{h\in B_s(R)}\Vert P_2^{j^\ast}f-h\Vert_{L_2}^2\leq \Vert P_2^{j^\ast}f-\widetilde{h}\Vert_{L_2}^2=:d^2,$$
		where
		\begin{align}
		d^2&= \Vert P_2^{j^\ast}f-\widetilde{h}\Vert_{L_2}^2\\
		&=\sum_{k=1}^{2^{\widetilde{j}}}\left(1-\frac{\sqrt{R^2-\Vert P_{j\neq\widetilde{j}}f\Vert_{\mathcal{B}_s}^2}}{\Vert P_{\widetilde{j}}f\Vert_{\mathcal{B}_s}}\right)^2a_{\widetilde{j},k}^2\\
		&=\left(1-\frac{\sqrt{R^2-\Vert P_{j\neq\widetilde{j}}f\Vert_{\mathcal{B}_s}^2}}{\Vert P_{\widetilde{j}}f\Vert_{\mathcal{B}_s}}\right)^2\frac{\Vert P_{\widetilde{j}}f\Vert_{\mathcal{B}_s}^2}{4^{\widetilde{j} s}}\\
		&=\left(\Vert P_{\widetilde{j}}f\Vert_{\mathcal{B}_s}-\sqrt{R^2-\Vert P_{j\neq\widetilde{j}}f\Vert_{\mathcal{B}_s}^2}\right)^24^{-\widetilde{j} s}.
		\end{align}
		This tells us that if $\Vert P_{j\neq\widetilde{j}}f\Vert_{\mathcal{B}_s}\leq R$,
		\begin{align}
		\Vert P_{\widetilde{j}}f\Vert_{\mathcal{B}_s}&=2^{\widetilde{j}s}d+\sqrt{R^2-\Vert P_{j\neq \widetilde{j}}f\Vert_{\mathcal{B}_s}^2}\geq 2^{\widetilde{j}s}d,\label{s:LBd1}\\
		\Vert P_2^{j^\ast}f\Vert_{\mathcal{B}_s}^2&=R^2+2\cdot 2^{\widetilde{j}s}d\sqrt{R^2-\Vert P_{j\neq\widetilde{j}}f\Vert_{\mathcal{B}_s}^2}+4^{\widetilde{j}s}d^2.\label{s:LBd2}\end{align}
		\item \underline{Bound in terms of $4^{j^\ast s}\rho_{j^\ast}^2$}\\
		If $\Vert P_2^{j^\ast-1}f\Vert_{\mathcal{B}_s}\leq R$, we can use \eqref{s:LBd2} with $\widetilde{j}=j^\ast$ and $d\geq\rho_{j^\ast}\geq 0$ and obtain
		$$\Vert P_2^{j^\ast}f\Vert_{\mathcal{B}_s}^2\geq R^2+4^{j^\ast s}\rho_{j^\ast}^2.$$
		If $\Vert P_2^{j^\ast-1}f\Vert_{\mathcal{B}_s}>R$, observe that by the triangle inequality
		$$\inf_{h\in B_s(R)}\Vert P_2^{j^\ast}f-h\Vert_{L_2}\leq \inf_{h\in B_s(R)}\Vert P_2^{j^\ast-1}f-h\Vert_{L_2}+\Vert P_{j^\ast}f\Vert_{L_2}\leq \rho_{j^\ast-1}+\Vert P_{j^\ast}f\Vert_{L_2}$$
		and since
		$$\rho_{j^\ast}-\rho_{j^\ast-1}\geq a\frac{2^{(3j^\ast+2J)/20}}{\sqrt{n}}(1-2^{-3/20})\geq \frac{1}{11}\rho_{j^\ast}=\frac{1}{A}\rho_{j^\ast},$$
		we obtain
		\begin{align}\Vert P_2^{j^\ast}f\Vert_{\mathcal{B}_s}^2&=\Vert P_2^{j^\ast-1}f\Vert_{\mathcal{B}_s}^2+\Vert P_{j^\ast}f\Vert_{\mathcal{B}_s}^2\\
		&\geq R^2+4^{j^\ast s}(\rho_{j^\ast}-\rho_{j^\ast-1})^2\\
		&\geq R^2+\frac{1}{A^2}4^{j^\ast s}\rho_{j^\ast}^2.\end{align}
		So, in any case,
		\begin{equation}
		\Vert P_2^{j^\ast}f\Vert_{\mathcal{B}_s}^2\geq R^2+\frac{1}{A^2}4^{j^\ast s}\rho_{j^\ast}^2.
		\label{s:simLB}
		\end{equation}
		\item \underline{Main arguments}~\\
		We are now ready to prove \eqref{s:LBH1} effectively. To that end, fix an index
		$$\overline{j}\in\underset{j\in\mathcal{J}}{\argmax}~2^{js}\Vert P_{j}f\Vert_{\mathcal{B}_s}.$$
		\begin{description}
			\item[\normalfont \textbf{Case 1:} $\boldsymbol{\Vert P_{j\neq\overline{j}}f\Vert_{\mathcal{B}_s}\leq R}$]~\\
			In that case, we can use \eqref{s:LBd1} and \eqref{s:LBd2} with $\widetilde{j}=\overline{j}$ in comination with \eqref{s:simLB} and obtain
			\begin{align}\Vert P_2^{j^\ast}f\Vert_{\mathcal{B}_s}^2&\geq R^2+2^{\overline{j}s}d\sqrt{R^2-\Vert P_{j\neq\overline{j}}f\Vert_{\mathcal{B}_s}^2}+\frac{1}{2}4^{\overline{j} s}d^2+\frac{1}{2\cdot A^2}4^{j^\ast s}d^2\\
			&=R^2+2^{\overline{j}s}d\left(\Vert P_{\overline{j}}f\Vert_{\mathcal{B}_s}-2^{\overline{j} s}d\right)+\frac{1}{2}4^{\overline{j} s}d^2+\frac{1}{2\cdot A^2}4^{j^\ast s}d^2\\
			&=R^2+d\cdot 2^{\overline{j}s}\left(\Vert P_{\overline{j}}f\Vert_{\mathcal{B}_s}-\frac{1}{2}2^{\overline{j} s}d\right)+\frac{1}{2\cdot A^2}4^{j^\ast s}d^2\\
			&\geq R^2+\rho_{j^\ast}2^{\overline{j}s}\left(\Vert P_{\overline{j}}f\Vert_{\mathcal{B}_s}-\frac{1}{2}2^{\overline{j} s}d\right)+\frac{1}{2\cdot A^2}4^{j^\ast s}\rho_{j^\ast}^2\\
			&\geq R^2+\frac{1}{2}\rho_{j^\ast}M_{j^\ast}+\frac{1}{2\cdot A^2}4^{j^\ast s}\rho_{j^\ast}^2,\end{align}
			remembering \eqref{s:MaxDef}.
			\item[\normalfont \textbf{Case 2:} $\boldsymbol{\Vert P_{j\neq\overline{j}}f\Vert_{\mathcal{B}_s}> R}$]~\\
			That case can be handled quickly by considering two subcases:
			\begin{description}
				\item[\normalfont \textbf{Subcase 1:} $\boldsymbol{4^{j^\ast s}\rho_{j^\ast}^2\geq \rho_{j^\ast}M_{j^\ast}}$]~\\
				Observe that with \eqref{s:simLB}
				\begin{align}\Vert P_2^{j^\ast}f\Vert_{\mathcal{B}_s}^2&\geq R^2+\frac{1}{2\cdot A^2}4^{j^\ast s}\rho_{j^\ast}^2+\frac{1}{2\cdot A^2}4^{j^\ast s}\rho_{j^\ast}^2\\
				&\geq R^2+\frac{1}{2\cdot A^2}\rho_{j^\ast}M_{j^\ast}+\frac{1}{2\cdot A^2}4^{j^\ast s}\rho_{j^\ast}^2.\end{align}
				\item[\normalfont \textbf{Subcase 2:} $\boldsymbol{4^{j^\ast s}\rho_{j^\ast}^2< \rho_{j^\ast}M_{j^\ast}}$]~\\
				In that case we have
				$$\Vert P_{\overline{j}}f\Vert_{\mathcal{B}_s}>\frac{4^{j^\ast s}}{2^{\overline{j}s}}\rho_{j^\ast}\geq 2^{\overline{j}s}\rho_{j^\ast}$$
				and thus
				\begin{align}
				\Vert P_2^{j^\ast}f\Vert_{\mathcal{B}_s}^2&=\Vert P_{j\neq\overline{j}}\Vert_{\mathcal{B}_s}^2+\Vert P_{\overline{j}}f\Vert_{\mathcal{B}_s}^2\\
				&>R^2+2^{\overline{j}s}\rho_{j^\ast}\Vert P_{\overline{j}}f\Vert_{\mathcal{B}_s}\\
				&= R^2+\rho_{j^\ast}M_{j^\ast}\\
				&\geq R^2+\frac{1}{2}\rho_{j^\ast}M_{j^\ast}+\frac{1}{2}4^{j^\ast s}\rho_{j^\ast}^2.
				\end{align}
			\end{description}
		\end{description}
		This concludes the proof since in any case \eqref{s:LBH1} holds.
	\end{enumerate}
\end{Proof}~\\
\subsection{Estimation of $\boldsymbol{M_{j^\ast}}$}\label{s:estimation}
As a last major step before directly controlling the type-I and type-II error probabilities, we need to find an appropriate estimator for $M_{j^\ast}$.
\begin{lemma}For $\delta\in(0,1)$ and $j^\ast,j\in\mathcal{J}$, let
	\begin{align}
	C_\delta&:=\sqrt{\frac{2}{\delta}},\\
	D_{j^\ast,\delta}&:=\frac{4^{j^\ast s}}{\sqrt{n}}\left(\sqrt{2}C_{\delta}+2^{j^\ast/4}\sqrt{C_{\delta}}\right)\\
	Y_j&:=16^{js}\left(\Vert P_j\widehat{f}\Vert_{L_2}^2-\frac{2^j}{n}\right)
	\end{align}
	and define the events
	\begin{align}\xi_{j^\ast,\delta}^0&:=\left\{M_{j^\ast}\leq \sqrt{\max_{2\leq j\leq j^\ast}|Y_{j}|}+D_{\delta,j^\ast}\right\},\label{s:Xi0}\\
	\xi_{j^\ast,\delta}^1&:=\left\{M_{j^\ast}\geq\sqrt{\max_{2\leq j\leq j^\ast}|Y_{j}|}-D_{\delta,j^\ast}\right\}.\label{s:Xi1}
	\end{align}	
	Then, for any monotone decreasing sequence $(\beta_j)_{j\in\mathcal{J}}$ in $(0,1)$, the following holds:
	\begin{equation}\PR(\xi_{j^\ast,\beta_{j^\ast}}^1)\geq 1-\sum_{j=2}^{j^\ast}\beta_{j},~~\PR(\xi_{j^\ast,\beta_{j^\ast}}^0)\geq 1-\beta_{j^\ast}.\label{s:resXi}\end{equation}
\end{lemma}
\begin{Proof}
	Remembering \eqref{s:eq:Chi2Para}, we know that for $j\in\{2,3,\ldots,j^\ast\}$
	$$Z_j:=4^{js}\Vert P_j\widehat{f}\Vert_{\mathcal{B}_s}^2=16^{js}\Vert P_j\widehat{f}\Vert_{L_2}^2$$
	has the properties
	\begin{align}
	\EW[Z_j]&=16^{js}\frac{2^j}{n}+4^{js}\Vert P_jf\Vert_{\mathcal{B}_s}^2,\\
	\Var[Z_j]&=2\cdot 16^{2js}\left(\frac{2^j}{n^2}+\frac{2}{n}\Vert P_jf\Vert_{L_2}^2\right)\\
	&=16^{js}\left(2\cdot 16^{js}\frac{2^j}{n^2}+\frac{4}{n}4^{js}\Vert P_jf\Vert_{\mathcal{B}_s}^2\right).
	\end{align}
	Now observe that for $\delta\in(0,1)$
	\begin{align}\sqrt{\frac{1}{\delta}\Var[Z_j]}&\leq\sqrt{\frac{2}{\delta}}2^{j/2}\frac{16^{js}}{n}+\frac{2}{\sqrt{\delta n}}4^{js}2^{js}\Vert P_jf\Vert_{\mathcal{B}_s}\\
	&\leq C_{\delta}2^{j^{\ast}/2}\frac{16^{j^\ast s}}{n}+\sqrt{2}C_{\delta}\frac{4^{j^\ast s}}{\sqrt{n}}M_{j^\ast}\\
	&=:v_{\delta,j^\ast}.
	\end{align}
	With $Y_j=Z_j-16^{js}\frac{2^j}{n}$, Chebyshev's inequality now tells us that
	\begin{equation}\PR\left(\left|Y_j-m_j^2\right|\geq v_{\delta,j^\ast}\right)\leq \delta.\label{s:EstimationIneq}\end{equation}
	We derive two bounds from this statement by lower bounding the the left hand side in two different ways:\\
	On the one hand, observe
	$$|Y_j-m_j^2|\geq ||Y_j|-m_j^2|\geq |Y_j|-m_j^2\geq |Y_j|-M_{j^\ast}^2.$$ 
	Now, since $(\beta_{j})_{j\in\mathcal{J}}$ is monotone decreasing, the sequence $(v_{\beta_j,j^\ast})_{j\in\mathcal{J}}$ is increasing, so that via a union bound we obtain
	\begin{align}
	\sum_{j=2}^{j^\ast}\beta_j &\geq \PR\left(\exists j\in\{2,3,\ldots,j^\ast\}:~|Y_j|\geq M_{j^\ast}^2+v_{\beta_j,j^\ast}\right)\\
	&\geq \PR\left(\exists j\in\{2,3,\ldots,j^\ast\}:~|Y_j|\geq M_{j^\ast}^2+v_{\beta_{j^\ast},j^\ast}\right)\\
	&=\PR\left(\sqrt{\max_{2\leq j\leq j^\ast}|Y_j|}\geq \sqrt{M_{j^\ast}^2+v_{\beta_{j^\ast},j^\ast}}\right).
	\end{align}
	With
	\begin{align}
	\sqrt{M_{j^\ast}^2+v_{\beta_{j^\ast},j^\ast}}&=\sqrt{\left(M_{j^\ast}+\frac{C_{\beta_{j^\ast}}}{\sqrt{2}}\frac{4^{j^\ast s}}{\sqrt{n}}\right)^2-\frac{C_{\beta_{j^\ast}}^2}{2}\frac{16^{j^\ast s}}{n}+C_{\beta_{j^\ast}}2^{j^{\ast}/2}\frac{16^{j^\ast s}}{n}}\\
	&\leq\sqrt{\left(M_{j^\ast}+\frac{C_{\beta_{j^\ast}}}{\sqrt{2}}\frac{4^{j^\ast s}}{\sqrt{n}}\right)^2+\frac{C_{\beta_{j^\ast}}^2}{2}\frac{16^{j^\ast s}}{n}+C_{\beta_{j^\ast}}2^{j^{\ast}/2}\frac{16^{j^\ast s}}{n}}\\
	&\leq M_{j^\ast}+\frac{4^{j^\ast s}}{\sqrt{n}}\left(\sqrt{2}C_{\beta_{j^\ast}}+2^{j^\ast/4}\sqrt{C_{\beta_{j^\ast}}}\right),
	\end{align}
	we have
	$$\PR\left(\sqrt{\max_{2\leq j\leq j^\ast}|Y_j|}\geq M_{j^\ast}+\frac{4^{j^\ast s}}{\sqrt{n}}\left(\sqrt{2}C_{\beta_{j^\ast}}+2^{j^\ast/4}\sqrt{C_{\beta_{j^\ast}}}\right)\right)\leq \sum_{j=2}^{j^\ast}\beta_j$$
	and hence the first claim from $\eqref{s:resXi}$.\\
	On the other hand, observe $$|Y_j-m_j^2|\geq m_j^2-|Y_j|$$
	and consider the specific case $j=\overline{j}$ in \eqref{s:EstimationIneq}:
	\begin{align}
	\beta_{j^\ast}&\geq \PR\left(|Y_{\overline{j}}|\leq M_{j^\ast}^2-v_{\beta_{j^\ast},j^\ast}\right)\\
	&\geq \PR\left(\max_{2\leq j\leq j^\ast}|Y_{j}|\leq M_{j^\ast}^2-v_{\beta_{j^\ast},j^\ast}\right)\\
	&=\PR\left(\max_{2\leq j\leq j^\ast}|Y_{j}|+\frac{16^{j^\ast s}}{n}\left(\frac{C_{\beta_{j^\ast}}^2}{2}+2^{j^\ast/2}C_{\beta_{j^\ast}} \right)\leq\left(M_{j^\ast}-\frac{C_{\beta_{j^\ast}}}{\sqrt{2}}\frac{4^{j^\ast s}}{\sqrt{n}}\right)^2\right)\\
	&\geq\PR\left(\sqrt{\max_{2\leq j\leq j^\ast}|Y_{j}|}+\frac{4^{j^\ast s}}{\sqrt{n}}\left(\sqrt{2}C_{\beta_{j^\ast}} +2^{j^\ast/4}\sqrt{C_{\beta_{j^\ast}}}\right)\leq M_{j^\ast}\right),
	\end{align}
	which asserts the second claim from \eqref{s:resXi}.
\end{Proof}
\subsection{Conclusion}\label{s:conclusion}
We will now assemble the individual results of the previous sections to obtain the claim of Theorem \ref{s:UB}. For $j\in\mathcal{J}$ we introduce
\begin{equation}
\rho_j=\frac{1346}{\sqrt{\eta}}\cdot \frac{2^{(3j+2J)/20}}{\sqrt{n}},~~\alpha_j=\eta\frac{1-2^{-1/5}}{4}2^{(j-J)/5},~~\beta_j=\eta\frac{1-2^{-1/2}}{2}2^{-j/2},\label{s:paraDefSpec}
\end{equation}
so that in particular 
$$\sum_{j=2}^J \alpha_j\leq\frac{\eta}{4},~~\sum_{j=2}^J \beta_j\leq\frac{\eta}{4}$$
and $(\beta_j)_{j\in\mathcal{J}}$ is monotone decreasing.
\paragraph{Result for fixed index}~\\
For $j^\ast\in\mathcal{J}$ define
$$T_{j^\ast,\alpha_{j^\ast}}=\Vert P_2^{j^\ast}\widehat{f}\Vert_{\mathcal{B}_s}^2-A_{j^\ast}-\frac{2}{\sqrt{\alpha_{j^\ast}}}\cdot \frac{\sqrt{j^\ast-1}}{\sqrt{n}}\sqrt{\max_{2\leq j\leq j^\ast}\left|Y_j\right|}.$$
Then under $H_0'\cap\xi_{j^\ast,\beta_{j^\ast}}^0$, \eqref{s:Xi0} and \eqref{s:eq:ConcIn2} yield that with probability at least $1-\alpha_{j^\ast}$
\begin{align}T_{j^\ast,\alpha_{j^\ast}}&\leq\Vert P_2^{j^\ast}f\Vert_{\mathcal{B}_s}^2+\frac{2}{\sqrt{\alpha_{j^\ast}}}\cdot \frac{\sqrt{j^\ast-1}}{\sqrt{n}}M_{j^\ast}+\frac{2}{\sqrt{\alpha_{j^\ast}}}4^{j^\ast s}\frac{2^{j^\ast/2}}{n}\\
&~~~~~~~~~~~~~~~~~~~~~~~~~~~~~~~~~~~~~~~~~~~~~~-\frac{2}{\sqrt{\alpha_{j^\ast}}}\cdot \frac{\sqrt{j^\ast-1}}{\sqrt{n}}\sqrt{\max_{2\leq j\leq j^\ast}\left|Y_j\right|}\\
&\leq R^2+\frac{2}{\sqrt{\alpha_{j^\ast}}}\cdot \frac{\sqrt{j^\ast-1}}{\sqrt{n}}D_{j^\ast,\beta_{j^\ast}}+\frac{2}{\sqrt{\alpha_{j^\ast}}}4^{j^\ast s}\frac{2^{j^\ast/2}}{n}
\end{align}
so that with
\begin{equation}\tau_{j^\ast,\alpha_{j^\ast}}=R^2+\frac{2}{\sqrt{\alpha_{j^\ast}}}\left( \frac{\sqrt{j^\ast-1}}{\sqrt{n}}D_{j^\ast,\beta_{j^\ast}}+4^{j^\ast s}\frac{2^{j^\ast/2}}{n}\right),\label{s:tauDef}\end{equation}
we obtain
\begin{equation}\PR_{H_0'}(T_{j^\ast,\alpha_{j^\ast}}>\tau_{j^\ast,\alpha_{j^\ast}}~|~\xi_{j^\ast,\beta_{j^\ast}}^0)\leq\alpha_{j^\ast}.\label{s:T1err}\end{equation}
On the other hand, let $\boldsymbol{j^{\ast}}$ be a transition index with property \eqref{s:Logic2}. Then under $H_1'\cap\xi_{\boldsymbol{j^{\ast}},\beta_{\boldsymbol{j^{\ast}}}}^1$, \eqref{s:eq:ConcIn2} and \eqref{s:LBH1} tell us that with probability at least $1-\alpha_{\boldsymbol{j^\ast}}$
\begin{align}T_{\boldsymbol{j^\ast},\alpha_{\boldsymbol{j^{\ast}}}}&\geq\Vert P_2^{\boldsymbol{j^\ast}}f\Vert_{\mathcal{B}_s}^2-\frac{2}{\sqrt{\alpha_{\boldsymbol{j^{\ast}}}}}\cdot \frac{\sqrt{\boldsymbol{j^{\ast}}-1}}{\sqrt{n}}M_{\boldsymbol{j^{\ast}}}-\frac{2}{\sqrt{\alpha_{\boldsymbol{j^{\ast}}}}}4^{\boldsymbol{j^{\ast}} s}\frac{2^{\boldsymbol{j^{\ast}}/2}}{n}\\
&~~~~~~~~~~~~~~~~~~~~~~~~~~-\frac{2}{\sqrt{\alpha_{\boldsymbol{j^{\ast}}}}}\cdot \frac{\sqrt{\boldsymbol{j^{\ast}}-1}}{\sqrt{n}}\sqrt{\max_{2\leq j\leq j^\ast}\left|Y_j\right|}\\
&\geq R^2+\left(\frac{1}{2\cdot \deno^2}\rho_{\boldsymbol{j^{\ast}}}-\frac{2}{\sqrt{\alpha_{\boldsymbol{j^{\ast}}}}}\cdot \frac{\sqrt{\boldsymbol{j^{\ast}}-1}}{\sqrt{n}}\right)M_{\boldsymbol{j^{\ast}}}+\frac{1}{2\cdot\deno^2}4^{\boldsymbol{j^{\ast}} s}\rho_{\boldsymbol{j^{\ast}}}^2\\
&~~~~~~~~~~~~~~~~~~~~~~~~~~-\frac{2}{\sqrt{\alpha_{\boldsymbol{j^{\ast}}}}}4^{\boldsymbol{j^{\ast}} s}\frac{2^{\boldsymbol{j^{\ast}}/2}}{n}-\frac{2}{\sqrt{\alpha_{\boldsymbol{j^{\ast}}}}}\cdot \frac{\sqrt{\boldsymbol{j^{\ast}}-1}}{\sqrt{n}}\sqrt{\max_{2\leq j\leq j^\ast}\left|Y_j\right|}.
\end{align}
Provided that
\begin{equation}\frac{1}{2\cdot \deno^2}\rho_{\boldsymbol{j^{\ast}}}\geq\frac{4}{\sqrt{\alpha_{\boldsymbol{j^{\ast}}}}}\cdot \frac{\sqrt{\boldsymbol{j^{\ast}}-1}}{\sqrt{n}}\label{s:TIIcond1},\end{equation}
using \eqref{s:Xi1} this yields
\begin{equation}
T_{\boldsymbol{j^\ast},\alpha_{\boldsymbol{j^{\ast}}}}\geq R^2+\frac{1}{2\cdot\deno^2}4^{\boldsymbol{j^{\ast}} s}\rho_{\boldsymbol{j^{\ast}}}^2-\frac{2}{\sqrt{\alpha_{\boldsymbol{j^{\ast}}}}}\cdot \frac{\sqrt{\boldsymbol{j^{\ast}}-1}}{\sqrt{n}}D_{\boldsymbol{j^{\ast}},\beta_{\boldsymbol{j^{\ast}}}}-\frac{2}{\sqrt{\alpha_{\boldsymbol{j^{\ast}}}}}4^{\boldsymbol{j^{\ast}} s}\frac{2^{\boldsymbol{j^{\ast}}/2}}{n}.\label{s:TLB1}
\end{equation}
Now by explicit computation we see that the choices in \eqref{s:paraDefSpec} ensure \eqref{s:TIIcond1} as well as
$$\frac{1}{4\cdot A^2}4^{\boldsymbol{j^{\ast}} s}\rho_{\boldsymbol{j^{\ast}}}^2\geq \frac{4}{\sqrt{\alpha_{\boldsymbol{j^{\ast}}}}}\cdot \frac{\sqrt{\boldsymbol{j^{\ast}}-1}}{\sqrt{n}}D_{\boldsymbol{j^{\ast}},\beta_{\boldsymbol{j^{\ast}}}}~~\mathrm{and}~~\frac{1}{4\cdot A^2}4^{\boldsymbol{j^{\ast}} s}\rho_{\boldsymbol{j^{\ast}}}^2\geq \frac{4}{\sqrt{\alpha_{\boldsymbol{j^{\ast}}}}}4^{\boldsymbol{j^{\ast}} s}\frac{2^{\boldsymbol{j^{\ast}}/2}}{n},$$
so that \eqref{s:TLB1} can be continued as
$$T_{\boldsymbol{j^{\ast}},\alpha_{\boldsymbol{j^{\ast}}}}\geq R^2+\frac{2}{\sqrt{\alpha_{\boldsymbol{j^{\ast}}}}}\left( \frac{\sqrt{\boldsymbol{j^{\ast}}-1}}{\sqrt{n}}D_{\boldsymbol{j^{\ast}},\beta_{\boldsymbol{j^{\ast}}}}+4^{\boldsymbol{j^{\ast}} s}\frac{2^{\boldsymbol{j^{\ast}}/2}}{n}\right)=\tau_{\boldsymbol{j^{\ast}},\alpha_{\boldsymbol{j^{\ast}}}}$$
and hence, finally,
\begin{equation}\PR_{H_1'}(T_{\boldsymbol{j^{\ast}},\alpha_{\boldsymbol{j^{\ast}}}}\leq \tau_{\boldsymbol{j^{\ast}},\alpha_{\boldsymbol{j^{\ast}}}}~|~\xi_{\boldsymbol{j^{\ast}},\beta_{\boldsymbol{j^{\ast}}}}^1)\leq\alpha_{\boldsymbol{j^{\ast}}}.\label{s:T1err}\end{equation}
\paragraph{Generalisation to unknown $\boldsymbol{j^\ast}$}~\\
For our test
\begin{equation}\varphi(P_2^J \widehat{f})=1-\prod_{j^\ast=2}^{J}\mathds{1}_{\{T_{j^\ast,\alpha_{j^\ast}}\leq\tau_{j^\ast,\alpha_{j^\ast}}\}},\label{s:thetest}\end{equation}
we can conclude with \eqref{s:resXi} and \eqref{s:paraDefSpec} that on the one hand
\begin{align}
\PR_{H_0'}\left(\varphi=1\right)&\leq \sum_{j^\ast=2}^J\left( \PR_{H_0'}\left(T_{j^\ast,\alpha_{j^\ast}}>\tau_{j^\ast,\alpha_{j^\ast}}~|~\xi_{j^{\ast},\beta_{j^\ast}}^0\right)+(1-\PR(\xi_{j^{\ast},\beta_{j^\ast}}^0))\right)\\
&\leq \frac{\eta}{4}+\frac{\eta}{4}=\frac{\eta}{2}\label{s:H0err}
\end{align}
and on the other hand
\begin{align}
\PR_{H_1'}\left(\varphi=0\right)&\leq \PR_{H_1'}\left(\forall j^\ast\in\mathcal{J}:~T_{j^\ast,\alpha}\leq\tau_{j^\ast,\alpha}\right)\\
&\leq \PR_{H_1'}\left(T_{\boldsymbol{j^\ast},\alpha}\leq\tau_{\boldsymbol{j^\ast},\alpha_{\boldsymbol{j^\ast}}}~|~\xi_{\boldsymbol{j^\ast},\beta_{\boldsymbol{j^\ast}}}^1\right)+\left(1-\PR(\xi_{\boldsymbol{j^\ast},\beta_{\boldsymbol{j^\ast}}}^1)\right)\\
&\leq \alpha_{\boldsymbol{j^\ast}}+\sum_{j=2}^{\boldsymbol{j^\ast}}\beta_{j}\\
&\leq \frac{\eta}{4}+\frac{\eta}{4}=\frac{\eta}{2}.\label{s:H1err}
\end{align}
\paragraph{Specification of $\boldsymbol{J}$ and conclusion}~\\
We are now ready to return to \eqref{s:rhoDef}. Choose
\begin{equation}J:=\left\lfloor \frac{1}{2t+1/2}\frac{\ln(n)}{\ln(2)}\right\rfloor,\label{s:Jdef}\end{equation}
so that
\begin{equation}\frac{1}{2}n^{\frac{1}{2t+1/2}}\leq 2^J\leq n^{\frac{1}{2t+1/2}}.\label{s:2Jbound}\end{equation}
That yields
\begin{equation}2^{-Jt}\leq 2^t n^{-\frac{t}{2t+1/2}}\label{s:2Jtbound}\end{equation}
and, on the other hand,
$$\rho_J=\frac{1346}{\sqrt{\eta}}\frac{2^{J/4}}{\sqrt{n}}\leq \frac{1346}{\sqrt{\eta}}\cdot n^{-\frac{t}{2t+1/2}}.$$
Therefore, whenever we choose
$$\rho\geq \left(\frac{1346}{\sqrt{\eta}}+\frac{R}{1-2^{-t}}\right)n^{-\frac{t}{2t+1/2}},$$
indeed by \eqref{s:H0err} and \eqref{s:H1err}
$$\sup_{f\in B_s(R)}\PR_f(\varphi=1)+\sup_{f\in \widetilde{B}_{s,t}(R,\rho)}\PR_f(\varphi=0)\leq \frac{\eta}{2}+\frac{\eta}{2}=\eta.$$
\section{Proof of Theorem \ref{s:LB}}
\subsection{Description of the Strategy}\label{LBderivation}
According to \eqref{s:testprob}, given $\eta\in(0,1)$, we aim at finding $\rho>0$ such that for any test $\varphi$,
$$\sup_{f\in B_s(R)}\PR(\varphi=1)+\sup_{f\in \widetilde{B}_{s,t}(R,\rho)}\PR(\varphi=0)>\eta.$$
This can be achieved through a Bayesian-type approach, see e.g. \cite{Bar02}: Let $\nu_0,\nu_\rho$ be probability distributions (priors) such that $\mathrm{supp}(\nu_0)\subseteq B_s(R)$ and $\mathrm{supp}(\nu_\rho)\subseteq \widetilde{B}_{s,t}(R,\rho)$. Then we have
\begin{align}
	\sup_{f\in B_s(R)}\PR_{f}(\varphi=1)&+\sup_{f\in\widetilde{B}_{s,t}(R,\rho)}\PR_{f}(\varphi=0)\\
	&~~~~~~~~~~~\geq\PR_{f\sim\nu_0}(\varphi=1)+\PR_{f\sim\nu_\rho}(\varphi=0)\\
	&~~~~~~~~~~~\geq 1-\frac{1}{2}\Vert \PR_{f\sim\nu_\rho}-\PR_{f\sim\nu_0}\Vert_{\mathrm{TV}}\\
	&~~~~~~~~~~~\geq 1-\frac{1}{2}\left(\int\left(\frac{\dx \PR_{f\sim\nu_\rho}}{\dx \PR_{f\sim\nu_0}}(x)\right)^2\ \dx \PR_{f\sim\nu_0}(x)-1\right)^{\frac{1}{2}}.\end{align}
This tells us that if we find $\widetilde{\rho}>0$ such that
\begin{equation}\int\left(\frac{\dx \PR_{f\sim \nu_\rho}}{\dx \PR_{f\sim\nu_0}}\right)^2\ \dx \PR_{f\sim\nu_0}< 1+4(1-\eta)^2,\label{Chi2eq}\end{equation}
for any test $\varphi$ it holds that
$$\sup_{f\in B_s(R)}\PR(\varphi=1)+\sup_{f\in \widetilde{B}_{s,t}(R,\rho)}\PR(\varphi=0)>\eta$$
and hence
$$\rho^\ast \geq \widetilde{\rho}.$$
\subsection{Application to our Problem}
\paragraph{Priors}~\\
Since the upper bound does not depend on $s$ and we found the index $J$ from \eqref{s:Jdef} to be critical, we choose the following structurally simple priors: Let $\nu_0$ be the Dirac-$\delta$ distribution on $\{0\}$ (i.e. $f\equiv 0$) and $\nu_\rho$ be the uniform distribution on
$$\mathcal{A}_{\rho,v}:=\left\{\sum_{k=1}^{2^J}a_{J,k}\psi_{j,k}~|~a_{J,1},a_{J,2},\ldots,a_{J,2^J}\in\{v,-v\}\right\},$$
where $v>0$ needs further specification: On the one hand, it is necessary to ensure that each $f\in\mathcal{A}_{\rho,v}$ fulfils $\Vert f\Vert_{\mathcal{B}_t}\leq R$ - note that for any such $f$, $\Vert f\Vert_{L_2}=2^{J/2}v$, so that by construction that condition reads $$2^{J(t+1/2)}v\leq R.$$
This motivates the choice $v:=a_\eta\cdot R\cdot 2^{-J(t+1/2)}$ for some $a_\eta\in(0,1]$ specified later based on further restrictions. On the other hand, we require
\begin{equation}
\rho\leq\inf_{h\in B_s(R)}\Vert f-h\Vert_{L_2}.\label{s:reqForV}
\end{equation}
Since only the level $J$ is involved, this is in fact merely the minimum over the Euclidean ball with radius $R\cdot 2^{-Js}$ so that
$$\inf_{h\in B_s(R)}\Vert f-h\Vert_{L_2}=\max\left(0,2^{J/2}v-R\cdot 2^{-Js}\right).$$
Now, by explicit computation we see that if
$$n\geq \left(\frac{2^{1+s-t}}{a_\eta}\right)^{\frac{2t+1/2}{s-t}},$$
with our choice of $v$ we have
$$\max\left(0,2^{J/2}v-R\cdot 2^{-Js}\right)\geq \frac{1}{2}2^{J/2}v=a_\eta\frac{R}{2}2^{-Jt},$$
so that \eqref{s:reqForV} holds if
\begin{equation}\rho\leq a_\eta\frac{R}{2}2^{-Jt}.\label{s:rhostructure}\end{equation}
\paragraph{Statistical distance}~\\
Again, the central task in this proof is to compute the $\chi^2$-divergence between $\PR_{f\sim\nu_0}$ and $\PR_{f\sim\nu_\rho}$. By construction, $\PR_{f\sim\nu_0}$ corresponds to the $2^J$-fold product of Gaussian distributions with mean $0$ and variance $\frac{1}{n}$, so that for $x\in\mathbb{R}^{2^J}$
$$\dx\PR_{f\sim\nu_0}(x)=\sqrt{\frac{n}{2\pi}}^{2^J}\prod_{k=1}^{2^J}\exp\left(-\frac{n}{2}x_k^2\right).$$
On the other hand, $\PR_{f\sim\nu_\rho}$ corresponds to a uniform mixture of $2^{2^J}$ products of $2^J$ independent Gaussians with means of the form $\pm v$ and variance $\frac{1}{n}$.\\~\\
Let $\mathcal{S}:=\{1,-1\}^{2^J}$ and $R$ be uniformly distributed on $\mathcal{S}$ (i.e. the product of $2^J$ Rademacher variables). Then
\begin{align}
\dx\PR_{f\sim \nu_\rho}(x)&=\frac{1}{2^{2^J}}\sum_{\alpha\in\mathcal{S}}\sqrt{\frac{n}{2\pi}}^{2^J}\prod_{k=1}^{2^J}\exp\left(-\frac{n}{2}\sum_{k=1}^{2^J}(x_k-\alpha_kv)^2\right)\\
&=\sqrt{\frac{n}{2\pi}}^{2^J}\EW_R\left[\prod_{k=1}^{2^J}\exp\left(-\frac{n}{2}(x_k-R_kv)^2\right)\right]\end{align}
and furthermore, with an independent copy $R'$ of $R$,
\begin{align}(\dx\PR_{f\sim \nu_\rho}(x))^2
&=\left(\frac{n}{2\pi}\right)^{2^J}\EW_{R,R'}\left[\prod_{k=1}^{2^J}\exp\left(-\frac{n}{2}\left[(x_k-R_kv)^2+(x_k-R'_kv)^2\right]\right)\right]\\
&=\left(\frac{n}{2\pi}\right)^{2^J}\exp\left(-2^Jnv^2\right)\EW_{R,R'}\left[\prod_{k=1}^{2^J}\exp\left(-nx_k^2+nvx_k(R_k+R_k')\right)\right].\end{align}
The quotient we need to integrate in \eqref{Chi2eq} therefore reads
\begin{align}\frac{(\dx\PR_{f\sim \nu_\rho})^2}{\dx \PR_{f\sim\nu_0}}(x)
&=\sqrt{\frac{n}{2\pi}}^{2^J}\exp\left(-2^Jnv^2\right)\EW_{R,R'}\left[\prod_{k=1}^{2^J}\exp\left(-\frac{n}{2}x_k^2+nvx_k(R_k+R_k')\right)\right]\\
&=\sqrt{\frac{n}{2\pi}}^{2^J}\exp\left(-2^Jnv^2\right)\\
&~~~~~~~~~\cdot\EW_{R,R'}\left[\prod_{k=1}^{2^J}\exp\left(-\frac{n}{2}\left(x_k-v(R_k+R_k')\right)^2\right)\exp\left(nv^2(1+R_kR_k')\right)\right]\\
&=\sqrt{\frac{n}{2\pi}}^{2^J}\EW_{R,R'}\left[\prod_{k=1}^{2^J}\exp\left(nv^2R_kR_k'\right)\prod_{k=1}^{2^J}\exp\left(-\frac{n}{2}\left(x_k-v(R_k+R_k')\right)^2\right)\right]
.\end{align}
Since the product of independent Rademacher variables is itself a Rademacher variable, we obtain
\begin{align}
\int_{\mathbb{R}^{2^J}}\frac{(\dx\PR_{f\sim \nu_\rho})^2}{\dx \PR_{f\sim\nu_0}}(x)~\dx x&=\EW_{R,R'}\left[\prod_{k=1}^{2^J}\exp\left(nv^2R_kR_k'\right)\right]\\
&=\EW_{R}\left[\prod_{k=1}^{2^J}\exp\left(nv^2R_k\right)\right]\\
&=\prod_{k=1}^{2^J}\EW_{R_k}\left[\exp\left(nv^2R_k\right)\right]\\
&=\left(\cosh(nv^2)\right)^{2^J}\\
&\leq\exp\left(2^J\frac{n^2v^4}{2}\right).
\end{align}
\paragraph{Conclusion}~\\
Now, \eqref{Chi2eq} holds if
$$\exp\left(2^J\frac{n^2v^4}{2}\right)< 1+4(1-\eta)^2$$
which, by explicit computation, is fulfilled if
$$a_\eta\leq \frac{2^{J(t+1/4)}}{\sqrt{n}R}\sqrt[4]{\ln(1+4(1-\eta)^2)}.$$
Through \eqref{s:2Jbound} and \eqref{s:2Jtbound} we find that
$$\frac{2^{J(t+1/4)}}{\sqrt{n}}\geq \frac{2^{-t}}{16}$$
and obtain the stronger condition
$$a_\eta\leq \frac{\sqrt[4]{\ln(1+4(1-\eta)^2)}}{2^t 16R}.$$
In summary: Let
$$a_\eta=\min\left\{1,\frac{\sqrt{\ln(1+4(1-\eta)^2)}}{2^t 16R}\right\}.$$ If
$$n\geq\left\lceil \left(\frac{2^{1+s-t}}{a_\eta}\right)^{\frac{2t+1/2}{s-t}}\right\rceil,$$
the priors $\nu_0$ and $\nu_\rho$ meet all requirements and the lower bound
$$\rho^\ast\geq a_\eta\frac{R}{2}2^{-Jt}\geq a_\eta\frac{R}{2}n^{-\frac{t}{2t+1/2}}$$
is established, where we write $C_{\eta}:=\frac{R}{2}a_\eta$.
\bibliographystyle{acm}
\bibliography{Biblio}
%

\end{document}